\documentclass[12pt,fleqn,twoside]{article}
\usepackage[cp1251]{inputenc}
\usepackage[ukrainian]{babel}
\usepackage{amsmath,amsfonts,amssymb}
\usepackage{latexsym,hhline,graphics,epsfig,wrapfig}

\newcounter{theorem}
\newcommand{\theor}{\par\refstepcounter{theorem}%
{\bf Теорема \thetheorem }.\,\,}

\newcounter{lemma}
\newcommand{\be}{\begin{equation}}
\newcommand{\ee}{\end{equation}}

\begin{document}

\large \pagestyle{empty}

\noindent {\small \textbf{УДК}\ \ 517.55}

\bigskip

A.К. Бахтин, Г.П. Бахтина, И.В. Денега \vskip 4mm

A.K. Bakhtin, G.P. Bakhtina, I.V. Denega

\vskip 4mm
\textbf{Оценки произведения внутренних радиусов взаимно неналегающих областей в многомерных комплексных пространствах}
\vskip 4mm
\textbf{Estimates of product of inner radii of mutually non-overlapping domains in multidimensional complex spaces}

\begin{flushright}
\parbox{16.5cm}{\footnotesize We study extremal problem
on the product of power of generalized inner radii of
non-overlapping domains in $\mathbb{C}^{n}$.

\vskip 1mm

Роботу присвячено розв'язанню екстремальної задачі про добуток
степенів узагальнених внутрішніх радіусів неналягаючих областей в
$\mathbb{C}^{n}$.}
\end{flushright}

\medskip

\vskip 3mm

\thispagestyle{empty}

Целью данной работы является рассмотрение задачи о произведении
степеней обобщенных внутренних радиусов полицилиндрических
неналегающих областей с полюсами на лучевой системе точек. Эта
задача относится к разряду задач с так называемыми "свободными"
полюсами (см., например, \cite{BBZ}). Пространственные аналоги
ряда известных результатов о неналегающих областях на плоскости
были получены в работе \cite{DP}. Для этого  в \cite{DP} было
обобщено понятие внутреннего радиуса, а именно, введено понятие
гармонического радиуса пространственной области $B \subset R^{n}$
относительно некоторой внутренней точки. Пожалуй, работа \cite{DP}
является единственной работой, где удалось значительно
продвинуться в получении результатов о неналегающих областях для
пространственного случая. В тоже время для случая комплексной
плоскости задачи о неналегающих областях представляют достаточно
хорошо разработанное направление геометрической теории функций
комплексного переменного (см., например [1--15]).

В работе \cite{Bah} был получен метод, который позволил обобщить
некоторые результаты геометрической теории функций комплексного
переменного на многомерные комплексные пространства. В частности,
в этой же работе были предложены аналоги известных теорем теории
однолистых функций.

Пусть $\mathbb{N}$, $\mathbb{R}$
-- множества натуральных и вещественных чисел соответственно,
$\mathbb{C}$ -- плоскость комплексных чисел,
$\overline{\mathbb{C}}=\mathbb{C}\bigcup\{\infty\}$ -- ее
одноточечная компактификация, $\mathbb{R^{+}}=(0,\infty)$. По
определению $\mathbb{C}^{n}=\underbrace{(\mathbb{C}\times\mathbb{C}
\times\ldots \times\mathbb{C}}_{n-\textrm{раз }}),$ $n\in
\mathbb{N}$ (см., например, [4--6]).

$\overline{\mathbb{C}}^{n}=\underbrace{(\overline{\mathbb{C}}\times\overline{\mathbb{C}}
\times\ldots \times\overline{\mathbb{C}})}_{n-\textrm{раз }}$ --
компактификация пространства $\mathbb{C}^{n}$, далее так
называемое пространство теории функций (см., например, [4--6]).

Ясно, что $\mathbb{C}^{1}=\mathbb{C}$,
$\overline{\mathbb{C}}^{1}=\overline{\mathbb{C}}$.

Бесконечно удаленными точками $\overline{\mathbb{C}}^{n}$ являются те точки,
у которых хотя бы одна координата бесконечна. Множество всех
бесконечно удаленных точек имеет комплексную размерность $n-1$.

Топология в $\overline{\mathbb{C}}^{n}$ вводится как в декартовом
произведении топологических пространств. В этой топологии
$\overline{\mathbb{C}}^{n}$ компактно (см., например, [4--6]).

Область $\mathbb{B}=B_{1}\times B_{2}\times\ldots\times
B_{n}\subset\overline{\mathbb{C}}^{n}$, где каждая область
$B_{k}\subset\overline{\mathbb{C}}$, $k=\overline{1, n}$
называется полицилиндрической областью в
$\overline{\mathbb{C}}^{n}$ (см., например, [4]). Области $B_{k}$,
$k=\overline{1, n}$ назовем координатными областями выше указанной полицилиндрической области
$\mathbb{B}$. В дальнейшем всюду для краткости будем обозначать полицилиндрическую область $\mathbb{B}$ через ее координатные области следующим образом $\mathbb{B}=\{B_{k}\}_{k=1}^{n}$.

Обобщенным внутренним $p$-радиусом, $p\in\mathbb{N}$, $p\leqslant n$
полицилиндрической области $\mathbb{B}$ в точке $\mathbb{A}$
($\mathbb{A}\in\mathbb{B}$), будем называть величину
\begin{center}
$\mathbb{R}_{p}(\mathbb{B},\mathbb{A}):=\left[\prod\limits_{k=1}^p
r(B_{k},a_{k})\right]^{\frac{1}{p}},$ $p\in\mathbb{N}$,
$p\leqslant n$
\end{center}
где $r(B_{k},a_{k})$ -- внутренний радиус координатной
области $B_{k}$ в точке $a_{k}$. Если $p=n$, тогда обобщенный
внутренний $p$-радиус будем называть просто обобщенным внутренним
радиусом
$$\mathbb{R}(\mathbb{B},\mathbb{A})=\left[\prod\limits_{k=1}^n
r(B_{k},a_{k})\right]^{\frac{1}{n}}.$$

Пусть $m,n\in \mathbb{N}$, $m\geqslant2$. Систему точек
$\Delta_{m}:=\left\{a_{k} \in \mathbb{C}:\,
k=\overline{1,m}\right\}$ назовем \textbf{\emph{$m$-лучевой}},
если $|a_{k}|\in\mathbb{R^{+}}$ при $k=\overline{1,m}$,
$$0=\arg a_{1}<\arg a_{2}<\ldots<\arg a_{m}<2\pi.$$

Систему точек $\{\mathbb{A}_{k}\}$ ($\mathbb{A}_{k}=\{a_{p}^{(k)}\}\in \mathbb{C}^{n}$),
$k=\overline{1,m}$ назовем \textbf{лучевой}, если при каждом фиксированном
$p_{0}$ последовательность $\{a_{p_{0}}^{(k)}\}$, $k=\overline{1,m}$, является $m$-лучевой системой
точек, $p_{0}=\overline{1,n}$.

В данной работе мы будем рассматривать лучевые системы точек следующего вида:

$$\mathbb{A}_{1}=\{1,1,\ldots,1\},\quad(\text{то есть}\quad
a_{p}^{(1)}=1,\quad p=\overline{1,n}),$$
\begin{equation}\label{1a}\end{equation}
$$\arg a_{p}^{(k)}<\arg a_{p}^{(k+1)},\quad
k=\overline{1,m-1},\quad \arg a_{p}^{(m)}<2\pi,\quad p=\overline{1,n}.$$

Система $\{\mathbb{B}_{k}\}$ ($\mathbb{B}_{k}=\{B_{p}^{(k)}\}_{p=1}^{n}$, $k=\overline{1,m}$)
называется системой полицилиндрических неналегающих областей, если
при каждом фиксированном $p_{0}$, $p_{0}=\overline{1,n}$, система областей
$\{B_{p_{0}}^{(k)}\}$, $k=\overline{1,m}$ является системой
неналегающих областей на $\overline{\mathbb{C}}$.

На комплексной плоскости $\mathbb{C}$ для произвольной $m$-лучевой системы точек
$\Delta_{m}=\{a_{k}\}_{k=1}^{m}$ и
$\gamma\in\mathbb{R^{+}}\cup\{0\}$ полагаем
$$L^{(\gamma)}(\Delta_{m}):=\prod\limits_{k=1}^n\left[
\chi\left(\Bigl|\frac{a_k}{a_{k+1}}\Bigr|^\frac{1}{2\alpha_k}\right)\right]^{1-\frac{1}{2}\gamma\alpha_k^2}
\prod\limits_{k=1}^n|a_k|^{1+\frac{1}{4}\gamma(\alpha_k+\alpha_{k-1})},$$
где $\chi(t)=\frac{1}{2}(t+t^{-1})$, $\alpha_{k}:=\displaystyle\frac{1}{\pi}\arg
\displaystyle\frac{a_{k+1}}{a_{k}},$ $\alpha_{m+1}:=\alpha_{1},$
$k=\overline{1, m}.$

Тогда для произвольной лучевой системы точек вида (1) обозначим
$$\mathcal{\mathbb{L}}^{(\gamma)}\left(\{\mathbb{A}_m\}\right):=\left\{L^{(\gamma)}\left(\{a_{1}^{(k)}\}_{k=1}^{m}\right),\: L^{(\gamma)}\left(\{a_{2}^{(k)}\}_{k=1}^{m}\right),\ldots, L^{(\gamma)}\left(\{a_{p}^{(k)}\}_{k=1}^{m}\right)\right\}.$$

Рассмотрим функционал
$$J_{m}(\gamma)=\mathbb{R}^{\gamma}(\mathbb{B}_{0},
\mathbb{A}_{0})\prod\limits_{k=1}^m \mathbb{R}(\mathbb{B}_{k},
\mathbb{A}_{k}),$$ где $\gamma\in\mathbb{R^{+}}$,
$\mathbb{A}_{0}=(0,0,\ldots,0)$,
$\mathbb{A}_{k}\in\mathbb{B}_{k}\subset\overline{\mathbb{C}}^{n}$,
$k=\overline{0,m}$ и система $\{\mathbb{B}_{k}\}_{k=0}^{m}$ является системой взаимно неналегающих полицилиндрических областей в $\overline{\mathbb{C}}^{n}$.

Тогда имеет место следующая теорема.
\textbf{\theor} Пусть $m,n\in \mathbb{N}$, $m\geqslant 5$,
$\gamma\in (0, \sqrt[3]{m}\,]$, $\mathbb{A}_{0}=(0,0,\ldots,0)$. Тогда
для произвольной лучевой системы точек вида (1)
$\{\mathbb{A}_{k}\}=\{a_{p}^{(k)}\}_{k=1}^{m}\in
\mathbb{C}^{n}$ такой, что
$\mathbb{L}^{(\gamma)}(\mathbb{A}_{n})=(1,1,\ldots,1)=\textbf{1}$,
$\mathbb{L}^{(0)}(\mathbb{A}_{n})=\textbf{1}$ и любого набора
взаимно непересекающихся полицилиндрических областей
$\mathbb{B}_k$, $\mathbb{A}_{k}\in\mathbb{B}_{k}\subset\overline{\mathbb{C}}^{n}$ ($k=\overline{0,m}$), справедливо следующее неравенство
$$
\mathbb{R}^{\gamma}(\mathbb{B}_{0},\mathbb{A}_{0})\prod\limits_{k=1}^m
\mathbb{R}(\mathbb{B}_{k},\mathbb{A}_{k})\leqslant\left(\frac{4}{m}\right)^{m}
\frac{(\frac{4\gamma}{m^{2}})^{\frac{\gamma}{m}}}{(1-\frac{\gamma}{m^{2}})
^{m+\frac{\gamma}{m}}}\left(\frac{1-\frac{\sqrt{\gamma}}{m}}{1+\frac{\sqrt{\gamma}}{m}}\right)^{2\sqrt{\gamma}}.
$$

\textbf{Доказательство теоремы 1.} Сделаем следующее
преобразование
$$\mathbb{R}^{\gamma}(\mathbb{B}_{0},\mathbb{A}_{0})\prod\limits_{k=1}^m
\mathbb{R}(\mathbb{B}_{k},\mathbb{A}_{k})=$$

$$=\left[\prod\limits_{p=1}^n
r(B_{p}^{(0)},a_{p}^{(0)})\right]^{\frac{\gamma}{n}}\prod\limits_{k=1}^m\left[\prod\limits_{p=1}^n
r(B_{p}^{(k)},a_{p}^{(k)})\right]^{\frac{1}{n}}=$$

$$=\left[\prod\limits_{p=1}^n
\left[r^{\gamma}(B_{p}^{(0)},a_{p}^{(0)})\prod\limits_{k=1}^m
r(B_{p}^{(k)},a_{p}^{(k)})\right]\right]^{\frac{1}{n}}.$$

Тогда для фиксированного $p=\overline{1,n}$ области $B_{p}^{(k)}$,
$k=\overline{0,m}$, образуют систему неналегающих областей на
$\overline{\mathbb{C}}$. Поэтому следуя работе \cite{Den}, имеем

$$r^{\gamma}(B_{p}^{(0)},a_{p}^{(0)})\prod\limits_{k=1}^m
r(B_{p}^{(k)},a_{p}^{(k)})\leqslant
\left(\frac{4}{m}\right)^{m}\frac{\left(\frac{4\gamma}{m^{2}}\right)
^{\frac{\gamma}{m}}}{\left(1-\frac{\gamma}{m^{2}}\right)^
{m+\frac{\gamma}{m}}}\left(\frac{1-\frac{\sqrt{\gamma}}{m}}{1+\frac{\sqrt{\gamma}}{m}}\right)^{2\sqrt{\gamma}}.$$

Отсюда имеем
$$\mathbb{R}^{\gamma}(\mathbb{B}_{0},\mathbb{A}_{0})\prod\limits_{k=1}^m
\mathbb{R}(\mathbb{B}_{k},\mathbb{A}_{k})\leqslant\left[\prod\limits_{p=1}^n
\left(\frac{4}{m}\right)^{m}\frac{\left(\frac{4\gamma}{m^{2}}\right)^{\frac{\gamma}{m}}}
{\left(1-\frac{\gamma}{m^{2}}\right)^{m+\frac{\gamma}{m}}}
\left(\frac{1-\frac{\sqrt{\gamma}}{m}}{1+\frac{\sqrt{\gamma}}{m}}\right)
^{2\sqrt{\gamma}}\right]^{\frac{1}{n}}=$$

$$=\left(\frac{4}{m}\right)^{m}
\frac{\left(\frac{4\gamma}{m^{2}}\right)^{\frac{\gamma}{m}}}
{\left(1-\frac{\gamma}{m^{2}}\right)^{m+\frac{\gamma}{m}}}
\left(\frac{1-\frac{\sqrt{\gamma}}{m}}{1+\frac{\sqrt{\gamma}}{m}}\right)^{2\sqrt{\gamma}}.$$

Теорема 1 доказана.

\begin{center} {\bf ЛИТЕРАТУРА}
\end{center}

\end{document}